\documentclass[12pt]{article}
\usepackage{graphicx,psfrag,epsfig, color,float}
\usepackage{amssymb,amsmath,amscd,amsthm}
\usepackage{verbatim}
\usepackage{graphicx,psfrag,epsfig}

\usepackage{graphicx}

\usepackage{tikz}
\def\mysuper#1#2{\sbox0{\ensuremath{#1}}\copy0\raise\ht0\hbox{\ensuremath{\scriptstyle #2}}}
\usepackage{comment}

\usepackage{xfrac}
\usepackage[active]{srcltx}
\usepackage{mathtools}
\usepackage{relsize}
\usepackage{bbm}
\usepackage[normalem]{ulem}

\newtheorem{theorem}{Theorem}[section]

\begin{document}
\title{A Criterion for Convergence of Diffusion Processes on Graphs}

\author{
Leonid Koralov\footnote{Dept of Mathematics, University of Maryland,
College Park, MD 20742, koralov@umd.edu},
Xiangyi Tao\footnote{Dept of Mathematics, University of Maryland,
College Park, MD 20742, xtao12@umd.edu} 
}

\date{}
\maketitle

\begin{abstract}
\
In this paper, we consider diffusion processes on a graph. On the $i$-th edge, each process is governed by an operator of the form  $D_{v_i}D^+_{u_i}$, and the behavior at the vertices is determined by the gluing conditions. The paper characterizes the weak convergence of a family of such processes to a limiting process in terms of the behavior of the scale and speed functions $u_i^n$ and $v_i^n$ and the parameters of the gluing conditions. This may be viewed as a generalization of Freidlin and Wentzell's result on the weak convergence of a family of processes on a line.
\end{abstract}



\section{Introduction}

 Consider a family of diffusion processes $X^n_t$ and another process $X_t$ with the same state space. It is an important problem in PDEs and probability theory to find conditions on the coefficients of the processes that would guarantee the convergence of their distributions. In particular, the coefficients of $X^n_t$ need not converge point-wise (as in certain homogenization results) for the convergence of processes to hold (see, e.g., \cite{BLP}, \cite{JKO}). In the one-dimensional case, for processes defined on the entire real line, necessary and sufficient conditions have been obtained in \cite{FW}. 
 Here, we study the convergence of processes on graphs. Markov processes on graphs emerge naturally as a result of averaging of randomly perturbed Hamiltonian systems, diffusion processes with invariant manifolds, etc. (see  \cite{FSW}, \cite{F1}). 

One-dimensional Markov processes with continuous trajectories on a segment, with the property that any point can be reached from any other point, have been described by Feller (see \cite{IM}. \cite{Mandl}). A similar classification of processes on graphs can be found in  \cite{FW2}. In the current paper, we have a sequence of such processes, and we prove a criterion for their convergence to a limiting Markov process on the same graph. In particular, we see how the presence of nearly-trapping regions near the vertex for the pre-limiting processes can result in extra stickiness at the vertex for the limiting one. Our result is a natural generalization of \cite{FW} from processes on the line to processes on graphs. 

\section{Problem Setup and the Formulation of the Result}
We are interested in the convergence of diffusion processes on graphs. By specifying one interior point per edge and introducing stopping times corresponding to the processes successively visiting such points, the question of weak convergence can be easily reduced to considering processes on star graphs, and thus we will focus on such a case.

Let $\mathcal{G}$ be a star graph with $k$ edges, denoted by $E_{1}, E_{2}, \dots, E_{k}$, sharing a common vertex $O$. 
Each edge $E_{i}$ is identified with the half-line $[0, \infty)$ and is parameterized by a spatial coordinate $x_{i} \in [0, \infty)$, where the vertex corresponds to the point $x_{i}=0$ for each~$i$.
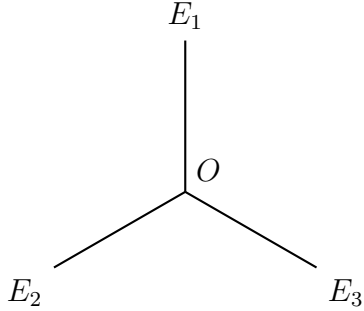
\begin{figure}[h!]
\centering
\begin{tikzpicture}[scale=2]
\draw[-, thick] (0,0) -- (0,1) node[above] {$E_{1}$};
\draw[-, thick] (0,0) -- (-0.87,-0.5) node[below left] {$E_{2}$};
\draw[-, thick] (0,0) -- (0.87,-0.5) node[below right] {$E_{3}$};

\node at (0,0) [above right] {$O$};
\end{tikzpicture}
\caption{A star graph with $k=3$ edges joined at the central vertex $O$.}
\label{fig:star_graph_n3}
\end{figure}

We will consider Markov processes on such graphs with the properties that all the trajectories are continuous, and each point can be reached from every other point. In particular, the vertex is accessible from each edge, and infinity is not reached in finite time. Each such process can be defined as follows. 
Consider strictly increasing functions $u_i, v_i$. $1 \leq i \leq k$, defined on $[0,\infty)$, such that

(a) The functions $u_i$ are continuous;

(b) The functions $v_i$ are right continuous;

(c) The functions satisfy $u_i(0) = v_i(0) = 0$ and $\int_0^\infty v_i d u_i = \infty$.
\\
Let $C_0(E_i)$ be the space of continuous functions on $E_i$ whose limit at infinity is equal to zero. The space of  continuous functions on $\mathcal{G}$ whose restrictions to $E_i$ belong to $C_0(E_i)$ will be denoted by $C_0(\mathcal{G})$. Suppose that there are positive numbers $\alpha_i$, $1 \leq i \leq k$,  and $\beta \geq 0$. For $x \in {\rm int}(E_i)$ and any function $f$, we define
\[
D_{u_i}^{+} f(x) = \lim_{h \downarrow 0} \frac{f(x+h) - f(x)}{{u_i}(x+h) - {u_i}(x)},\quad D_{v_i} f(x) = \lim_{h \to 0} \frac{f(x+h) - f(x)}{{v_i}(x+h) - {v_i}(x)}.
\]
Let $\mathcal{A}_i f = D_{v_i} D_{u_i}^+ f$ be the operator defined on such functions $f \in C_0(E_i)$ for which the derivatives can be successively taken in the interior of the edge and the resulting expression can be extended to an element of $C_0(E_i)$. The latter condition can be shown to imply that $D^+_{u_i} f$    and $D^+_{v_i} D^+_{u_i} f$ are defined on $E_i$ (including the origin) and the latter is continuous.

 Define the operator $\mathcal{A}$ on $C_0(\mathcal{G})$ as follows. The domain $D(\mathcal{A})$ consists of those functions $f \in C_0(\mathcal{G})$ for which $\mathcal{A}_i f$, $1 \leq i \leq k$,  can be viewed as an element of $C_0(\mathcal{G})$ (in particular, $\mathcal{A}_i f(0)$ does not depend on $i$ and can be denoted by $\mathcal{A}f(0)$), and the following gluing condition holds: 
 \begin{equation} \label{glco}
\sum_{i=1}^{k}\alpha_{i}D^+_{u_{i}} f(0)=\beta \mathcal{A} f (0).
\end{equation}
By Theorem 3.1 of \cite{FW2}, 
there exists a Markov process on
$\mathcal G$ with continuous sample paths whose generator is $\mathcal{A}$, and the distribution of such a process is determined uniquely. Under our assumption
that infinity is inaccessible on each edge, this process is
conservative. We denote such a process by $X_t$. We will consider not one such process but a sequence of processes $X^n_t$ with the generators $\mathcal{A}^n$, each defined as above, with the functions $u^n_i, v^n_i$, and the coefficients $\alpha^n_i$ and $\beta^n$, respectively. The goal is to identify the conditions on the functions and the coefficients in the gluing condition that would lead to the convergence of $X^n_t$ to $X_t$ (in the sense of the weak convergence of the induced measures on $C([0,\infty), \mathcal{G})$, for each initial point).

In the case where the processes $X^n_t$ and $X_t$ are defined on a line rather than a graph (or $X_t$ is defined on an open segment), the problem has been solved in \cite{FW}. The result does not involve gluing conditions, since there is no vertex in that case. Of course, a line can be viewed as a graph (with the origin serving as a vertex connecting two semi-infinite edges), and so, formally speaking, our result extends that of \cite{FW}. 

Processes on a segment (or half-line) with reflection or delay at one of the endpoints were considered in \cite{Kim}. A half-line is, in our terminology, just a graph with one edge. Thus, our results for the convergence of processes on graphs also apply to processes on a segment or half-line. Note that we describe an additional mechanism for convergence, compared to \cite{Kim}, where a delay at the vertex may appear in the limit, while the processes $X^n_t$ have $\beta_n = 0$ (see the example below our main result).

Observe that the functions $u_i, v_i$ and the coefficients in the gluing condition are, in general, not defined uniquely. Indeed, we could replace $\alpha_i$ by $a_i \alpha_i$, $u_i$ by $a_i u_i$, and $v_i$ by $v_i/a_i$, where $a_i$ is an arbitrary positive constant, resulting in the same operator.  Thus, without loss of generality, we can assume that $\alpha_i = 1/k$ for each $i$. We can also assume that $\alpha^n_i = 1/k$ for all $n$ and $i$ (it is not difficult to show that processes with $\alpha^n_i = 0$ cannot converge to a process with $\alpha_i \neq 0).$
Thus, the gluing conditions become
 \begin{equation}
\frac{1}{k} \sum_{i=1}^{k} D^+_{u_{i}} f(0)=\beta \mathcal{A} f (0). \label{glc1}
\end{equation}
 \begin{equation}
\frac{1}{k} \sum_{i=1}^{k} D^+_{u^n_{i}} f(0)=\beta^n \mathcal{A}^n f (0). \label{glc2}
\end{equation}

Our main result is the following theorem.
\begin{theorem} \label{mthe}
Suppose that the functions $u_i, v_i, u^n_i, v^n_i$, $n \geq 1$, $1 \leq i \leq k$, satisfy the assumptions stated above, and $\beta, \beta^n$ are nonnegative. Let $X_t$, $X^n_t$  be the corresponding Markov processes on $\mathcal{G}$ with the gluing conditions \eqref{glc1}, \eqref{glc2}, respectively. Then $X^n_t \rightarrow X_t$ (for each initial point) if and only if
there exist constants $c^{n}_{i}$ such that
\[
u^{n}_{i}(x)\longrightarrow u_{i}(x)\quad\text{for every }x \geq 0, \tag{C1}
\]
\[
v^{n}_{i}(x)-c^{n}_{i}\longrightarrow v_{i}(x)\quad\text{for every}~x > 0~\text{that is a continuity point of }v_{i}, \tag{C2}
\]
\[
\beta^{n}+\frac{1}{k}\sum_{i=1}^{k}c^{n}_{i}\longrightarrow \beta~~~as~n \rightarrow \infty. \tag{C3}
\]

\end{theorem}
\vspace{0.5cm}
\noindent
{\bf Example}.
Consider the one-edge special case $\mathcal G=[0,\infty)$ with the vertex $O=0$. Let
\[
u_1(x)=u_1^n(x)=x,\qquad v_1(x) = x,\qquad
v^n_1 = 
\begin{cases}
nx,~~~~~~~~~~~~~~0 \leq x \leq 1/n,\\
1+x-1/n,~~x \geq 1/n,\\
\end{cases} 
\]
while
\[
c^n_1 = 1, \qquad \beta = 1,\qquad \beta^n = 0.
\]
With this choice of functions and gluing conditions, the assumptions of Theorem~\ref{mthe} are satisfied. Note that the processes $X^n_t$ 
are instantaneously reflecting at $O$, while the limiting process has a delay at the origin.

\section{Proof of the Result}
Here, we prove Theorem~\ref{mthe}.
Let us prove the sufficiency part first. Let $\varphi_{i} \in E_i$, $\varphi_i > 0$.  We will treat $\varphi_i$ as a positive real number and as an element of $E_i$, depending on the context. Let $\Phi = \bigcup_i \{\varphi_i\}$ and $U = \bigcup_i \{x \in E_i, x \leq \varphi_i\}$. For the processes $X_t$ and $X^n_t$ starting at the origin, we will calculate the probability of reaching $\varphi_{i} \in E_i$ prior to reaching $\Phi \setminus \{\varphi_i\}$. We will also calculate the expectation of the time needed to reach $\Phi$.  

Let $p_i(x) = \mathrm{P}_x(X_{\sigma(\Phi)} = \varphi_i)$, $x \in U$,
where $\sigma(A)$ is the first time the process reaches a set $A$. 
Then $\mathcal{A}_j p_i(x) = 0$, $x \in (0, \varphi_j)$, for all $j$; $p_i(\varphi_i) =  1$; $p_i(\varphi_j) = 0$ for $j \neq i$; for each $i$,  $p_i$ satisfies the gluing condition \eqref{glco} at the vertex. This is justified as follows. Let $g_i$ solve the following system of equations 
\[
 \mathcal{A}_j g_i(x) = 0,~~x \in (0, \varphi_j)~~{\rm for~all}~~j; ~~~g_i(\varphi_i) =  1;~ g_i(\varphi_j) = 0~~{\rm  for }~ j \neq i;
 \]
 \[   \sum_{j=1}^{k}\alpha_{j}D^+_{u_{j}} g_i(0)=0
\]
 in the class of continuous functions on $U$. We extend $g_i$ to $\mathcal{G}$ as an element of $D(\mathcal{A})$. Apply Dynkin's formula to the function $g_i$ and the stopped process $X_{t \wedge \sigma(\Phi)}$ starting at $x \in U$:
 $$ \mathrm{E}_x[g_i(X_{t \wedge \sigma(\Phi)})] = g_i(x) + \mathrm{E}_x\left[\int_{0}^{t \wedge \sigma(\Phi)}  Ag_i(X_s) ds\right] = g_i(x), $$
 which gives $p_i(x) = g_i(x)$, $x \in U$, since the boundary data for $p_i$ and $g_i$ agree.

For $x \in [0, \varphi_i]$, let $q_i(x) = \mathrm{P}_x(X_{\sigma(\{0, \varphi_i\})} = \varphi_i) = u_i(x)/u_i(\varphi_i)$.  For $x \in E_j$ with $j \neq i$, $p_i(x) = p_i(0) (1 - q_j(x))$ by the strong Markov property, while $p_i(x) = p_i(0) + (1- p_i(0)) q_i(x)$ for $x \in E_i$. The gluing condition gives
\[
D_{u_i}p_i(0) + p_i(0) \sum_{j \neq i} D_{u_j} (1 - q_j) (0) = 0,
\]
and so
\[
\frac{1-p_i(0)}{u_i(\varphi_i)}- p_i(0) \sum_{j \neq i}  \frac{1}{u_j(\varphi_j)}= 0.
\]
Also, $\sum_i p_i(0) = 1$. Solving the system gives
\begin{equation*}
    p_i(0) = \frac{ 1/ u_i(\varphi_i) }{ \sum_{j=1}^k 1/ u_j(\varphi_j) }, ~~1 \leq i \leq k. 
\end{equation*}

Let $T(x) = \mathrm{E}_x \sigma(\Phi)$ and, for
 $x \in [0, \varphi_i]$, let $T_i(x) = \mathrm{E}_x \sigma(\{0,  \varphi_i\})$ be the first time the process reaches one of the endpoints of the segment $[0,\varphi_i]$. Then 
$\mathcal{A}_i T_i (x) = -1$,
$x \in (0, \varphi_i)$; $T_i(0)=T_i(\varphi_i) = 0$. Thus,
\[
D_{u_i} T_i(x) = -v_{i}(x_{i}) + \frac{1}{u_i(\varphi_i)} \int_0^{\varphi_{i}} v_{i}(y)du_{i}(y). 
\]

By the strong Markov property, $
T(x) = T_i(x) + (1-q_i(x)) T(0)$ for $x \in [0, \varphi_i]$. Similarly to the way done for $p_i$, it is easily shown that
$T(x)$ satisfies the gluing condition
\[
\frac{1}{k} \sum_{i=1}^{k} D^+_{u_{i}} T(0)+\beta =0,
\]
and so
\[
\frac{1}{k} \sum_{i=1}^{k} \frac{1}{u_i(\varphi_i)} \left( \int_0^{\varphi_{i}} v_{i}(y)du_{i}(y)  - {T(0)}\right)+\beta =0.
\]
Solving for $T(0)$ gives 
\begin{equation} 
\label{formulaforT} 
T(0) = \left( \sum_{i=1}^k \frac{1}{u_i(\varphi_i)} \int_0^{\varphi_i} v_i(y) du_i(y) + k\beta \right) \left( \sum_{i=1}^k \frac{1}{u_i(\varphi_i)} \right)^{-1}.
\end{equation}
We derived the formulas for the exit probabilities and average exit times for the process $X_t$. Similar formulas apply to the process $X^n_t$, and the corresponding quantities will be denoted by $p^n_i$ and $T^n$.

Let us fix $t \geq 0$ and $f \in D(\mathcal{A})$. To prove the convergence of the processes, it is sufficient  to show that
\begin{equation} \label{exp11}
\mathrm{E}_x\left(f(X^n_t)-f(X^n_0)-\int_{0}^{t}\mathcal{A}f(X^n_s)ds\right) \rightarrow 0~~~{\rm as}~~n \rightarrow \infty
\end{equation}
uniformly in $x \in K$ for any compact set $K$. This sufficiency is a version of Lemma 3.1, Chapter 8, in \cite{FW}. It requires the tightness of $X^n_s$, which is obvious here. We will fix $x$ in the arguments below, keeping in mind that all the estimates are uniform for $x$ in a compact. 

To demonstrate this convergence, let us fix $\varepsilon > 0$. Let $\eta > 0$ (to be selected later), and $\varphi_{i} = \varphi_i(\eta):= {u_{i}}^{-1}(\eta)$ for each $i$. We will split the interval $[0,t]$ into the intervals it takes the process $X^n_t$ to go from the origin to $\Phi = \Phi(\eta) = \bigcup_i \{\varphi_i(\eta)\}$ and back. By selecting a sufficiently small $\eta$, we will make the absolute value of the expression in the left-hand side of \eqref{exp11} smaller than $\varepsilon$ for all sufficiently large $n$.

We define the sequence of stopping times as follows: let $\tau^n_0 = \sigma^n_0 = 0$, and for $m \geq 1$, define recursively:
\begin{align*}
\tau^n_m &= \inf \{t \geq \sigma^n_{m-1} : X^n_t = 0 \}, \\
\sigma^n_m &= \inf \{t \geq  \tau^n_m : X^n_t \in \Phi \}.
\end{align*}
An interval $[\sigma^n_{m-1}, \tau^n_m]$ represents an inward excursion  (from the set $\Phi$ to the vertex $O$ or from the initial point to the vertex (if $m = 1$)).
An interval $[\tau^n_m, \sigma^n_m]$ represents an outward excursion (from the vertex $O$ to the set $\Phi$).

Let $M = M(t) = \max(m: \tau^n_m \leq t)$. Thus, $M = 0$ if the process doesn't reach $O$ prior to time $t$, and, otherwise, the random variable $M-1$ is equal to the number of completed outward-inward excursion pairs prior to time $t$.
%
%
Observe that $\tau^n_M$ is not a stopping time, but 
$\tau^n_{M+1}$ is. 
We rewrite the expectation in \eqref{exp11} as

\begin{align} 
    & \mathrm{E}_x\left(f(X^n_t)-f(X^n_0)-\int_{0}^{t}\mathcal{A}f(X^n_s)ds\right) \label{inward}\\
    =& \mathrm{E}_x\sum_{m=0}^{M} \left(f(X^n_{\tau^n_{m+1} \wedge t})-f(X^n_{\sigma^n_m \wedge t})-\int_{\sigma^n_m \wedge t}^{\tau^n_{m+1} \wedge t }\mathcal{A}f(X^n_s)ds\right) \label{1st-term} \\
    +& \mathrm{E}_x\sum_{m=1}^{M} \left( f(X^n_{\sigma^n_m})-f(X^n_{\tau^n_m})-\int_{\tau^n_m}^{\sigma^n_m}\mathcal{A}f(X^n_s)ds\right) \label{2nd-term}\\
    -& \mathrm{E}_x\left(\mathbf{1}_{\tau^n_M \leq t < \sigma^n_M  }\left(f(X^n_{\sigma^n_{M}})-f(X^n_t)-\int_{t}^{\sigma^n_{M}}\mathcal{A}f(X^n_s)ds\right) \right). \label{3rd-term}
\end{align}

We will estimate each of the three terms in the right-hand side separately. To bound \eqref{3rd-term}, we treat $X_t^n$ as a new starting point, which may be assumed to belong to $U$. Thus, the absolute value of  \eqref{3rd-term} is bounded from above by 
\begin{align}
    \sup_{x\in U} \left| \mathrm{E}_x\left(f(X^n_{\sigma^n})-f(x)-\int_{0}^{\sigma^n}\mathcal{A}f(X^n_s)ds\right) \right| \label{3rd bound},
\end{align}
where $\sigma^n$ is the first time the process reaches $\Phi$.
The expression 
$ \mathrm{E}_x\left|f(X^n_{\sigma^n})-f(x)\right|$
can be made arbitrarily small, for all $n$ and $x$, by selecting a sufficiently small $\eta$, due to the continuity of $f$.  We also have
\[
\left|
\mathrm{E}_x \int_{0}^{\sigma^n}\mathcal{A}f(X^n_s)ds \right| \leq \sup \mathcal{A} f \cdot
\mathrm{E}_x \sigma^n.
\]
The  expectation in the right-hand side is bounded by $\mathrm{E}_x \sigma^n (\{O\} \bigcup \Phi) + T^n(0)$. The term $\mathrm{E}_x \sigma^n (\{O\} \bigcup \Phi)$ can be easily expressed in terms of the functions $u^n_i$ and $v^n_i$, and can be made arbitrarily small by selecting sufficiently small $\eta$. It is shown below, when we prove \eqref{2nd-term}, that $T^n(0)$ can also be made arbitrarily small. This shows that \eqref{3rd bound} (and consequently  \eqref{3rd-term}) is bounded from above by $\varepsilon/5$ if $\eta$ is sufficiently small.

Next, let us deal with the second term, \eqref{2nd-term}, which is the most difficult one. It is equal to 
\begin{equation} \label{twofactors}
\mathrm{E}_x M \cdot \mathrm{E}_O \left(f(X^n_{\sigma^n})-f(X^n_0)-\int_{0}^{\sigma^n}\mathcal{A}f(X^n_s)ds\right).
\end{equation}
 Let us estimate $\mathrm{E}_x M$ from above. On each edge $E_i$, select a point $r_i > \varphi_i$ at a positive distance from the origin. For the process $X^n_t$ starting at $\varphi_i(\eta)$, we have 
\[ 
    \mathrm{P}_{\varphi_i}(X^n_{\sigma^n(\{0, r_i\})} =r_i) = u^n_i(\varphi_i)/u^n_i(r_i).
\]
Observe that there is $c > 0$ such that
\[
\mathrm{P}_{r_i} (\tau^n_1 \geq t) \geq c
\]
for each $i$ and $n$, and therefore
\[  
    \mathrm{P}_{\varphi_i} (\tau^n_1 \geq t) \geq c u^n_i(\varphi_i)/u^n_i(r_i).
\]
In other words, with a small probability (yet bounded from below), it takes at least time  $t$ for $X^n_t$ to make one excursion. 
This implies that 
\[ 
    \mathrm{P}_x (M > m) \leq \max_i  \left(1 - c \frac{u^n_i(\varphi_i)}{u^n_i(r_i)}\right)^m,~~~m \geq 1. 
\]
Therefore, for some $C > 0$,
\begin{equation}  
    \mathrm{E}_x M \leq C \max_i u^n_i (r_i)/ u^n_i (\varphi_i) \label{m_boundary}. 
\end{equation}

Now, let us bound the second factor in \eqref{twofactors}. We have
\begin{align}
    &\mathrm{E}_O \left(f(X^n_{\sigma^n})-f(X^n_0)-\int_{0}^{\sigma^n}\mathcal{A}f(X^n_s)ds\right) \notag\\
    =& \sum_{i=1}^{k} p^n_i \left(f(\varphi_i)-f(0)\right) - \mathcal{A}f(0) T^n(0) + \mathrm{E}_O    \int_{0}^{\sigma^n}( \mathcal{A}f(0) - \mathcal{A}f(X^n_s))ds. \label{abc12}
\end{align}

The absolute value of the last term on the right-hand side can be estimated from above by $T^n(0) a_1(\eta)$ with
\[
a_1(\eta) = \sup_{x \in \mathcal{G}_\eta} |\mathcal{A}f(0) - \mathcal{A}f(x)|,
\]
where $x \in \mathcal{G}_\eta$ if $x \in E_i$ for some $i$ and $0 \leq x \leq \varphi_i(\eta)$. Observe that $a_1(\eta) \rightarrow 0$ as $\eta \downarrow 0$ since $\mathcal{A}f$ is continuous. 

Since $v_j^n$ is increasing, the first term of the numerator in formula \eqref{formulaforT} for $T^n$ is bounded by $\sum_{j}v_j^n(\varphi_j)$, and therefore
$$T^n(0) \le \left({\sum_{j=1}^k v_j^n(\varphi_j) + k\beta^n}\right)\left({\sum_{j=1}^k \frac{1}{u_j^n(\varphi_j)}}\right)^{-1}.$$

Now, consider the product $\mathrm{E}_x M \cdot T^n(0)$. Using the bound for $\mathrm{E}_x M$:
\begin{align}
\mathrm{E}_x M \cdot T^n(0) &\le \left( C \max_i \frac{u_i^n(r_i)}{u_i^n(\varphi_i)} \right) \cdot \left({\sum_{j=1}^k v_j^n(\varphi_j) + k\beta^n}\right)\left({\sum_{j=1}^k \frac{1}{u_j^n(\varphi_j)}}\right)^{-1}\notag \\ \nonumber
&\le C \max_i u_i^n(r_i) \cdot \left( \sum_{j=1}^k v_j^n(\varphi_j) + k\beta^n \right)  .
\end{align}
Take a constant $C_1 > 0$ such that
\[
C \max_i u_i(r_i) \cdot \left( \sum_{j=1}^k v_j(r_j) + k\beta \right) \leq C_1
\]
for all $\eta \leq 1$. Therefore, $\mathrm{E}_x M \cdot T^n(0) \leq 2C_1$ for all sufficiently small $\eta$ and all sufficiently large $n$ (depending on $\eta$). Therefore,
\begin{equation}
    \left|\mathrm{E}_x M \cdot \mathrm{E}_O    \int_{0}^{\sigma^n}( \mathcal{A}f(0) - \mathcal{A}f(X^n_s))ds \right|\leq 2 C_1 a_1(\eta) \leq \frac{\varepsilon}{5} ,\label{eq:19}
\end{equation}
where the last inequality holds for all sufficiently small $\eta$, and the first inequality holds
for all sufficiently large $n$ (depending on $\eta$). 

Now, let us bound the product of $\mathrm{E}_x M$ and the first two terms on the right-hand side of \eqref{abc12}. Denote
\[
\bar v_j^n:=v_j^n-c_j^n.
\]
Then, since $u_j^n(0)=0$,
\[
\frac{1}{u_j^n(\varphi_j)}
\int_0^{\varphi_j} v_j^n(y)\,du_j^n(y)
=
\frac{1}{u_j^n(\varphi_j)}
\int_0^{\varphi_j} \bar v_j^n(y)\,du_j^n(y)
+c_j^n.
\]
Hence,
\begin{align*}
& \mathrm{E}_x M\cdot
\left|
\sum_{i=1}^{k} p_i^n(f(\varphi_i)-f(0))
-\mathcal{A}f(0)T^n(0)
\right|
\\
=&\,
\mathrm{E}_x M\cdot
\left|
\sum_{i=1}^{k}
\frac{1/u_i^n(\varphi_i)}
{\sum_{j=1}^k1/u_j^n(\varphi_j)}
(f(\varphi_i)-f(0))
\right.\\
&\left.
\qquad
-\mathcal{A}f(0)
\left(
\sum_{j=1}^k
\frac{1}{u_j^n(\varphi_j)}
\int_0^{\varphi_j}\bar v_j^n(y)\,du_j^n(y)
+\sum_{j=1}^k c_j^n+k\beta^n
\right)
\left(
\sum_{j=1}^k\frac{1}{u_j^n(\varphi_j)}
\right)^{-1}
\right|\\
=:&\,
\mathrm{E}_x M\cdot|\mathcal E^n|
\nonumber\\ 
\leq&
\mathrm{E}_x M\cdot
|\mathcal E^n-\mathcal E^\infty(\eta)|
+
\mathrm{E}_x M\cdot
|\mathcal E^\infty(\eta)|.
\end{align*}


where
\[
  \mathcal{E}^\infty(\eta) =  \frac{1}{k}\sum_{i=1}^k (f(\varphi_i)-f(0)) - \mathcal{A}f(0) \left( \frac{1}{k}\sum_{j=1}^k \int_0^{\varphi_j} v_j du_j + \beta\eta \right).
\]
Since $D^+_{u_i} f(0)$ is defined and by the gluing condition \eqref{glc1}, $\mathrm{E}_x M\cdot |\mathcal{E}^\infty(\eta)|$ is bounded by
\begin{align}
    \mathrm{E}_x M\cdot|\mathcal{E}^\infty(\eta)|  &\leq C \max_i \frac{u^n_i (r_i)}{u^n_i(\varphi_i)}  \left| \frac{1}{k}\sum_{i=1}^{k} (D_{u_i}^+ f(0)\eta + o(\eta)) - \mathcal{A}f(0) \left( \frac{1}{k}\sum_{j=1}^k \int_0^{\varphi_j} v_j du_j + \beta\eta \right) \right|\nonumber \\
    &\leq 2C \max_i \frac{u_i (r_i)}{u_i(\varphi_i)}  \left| \frac{1}{k}\sum_{i=1}^{k} (D_{u_i}^+ f(0)\eta + o(\eta)) - \mathcal{A}f(0) \left( \frac{1}{k}\sum_{j=1}^k \int_0^{\varphi_j} v_j du_j + \beta\eta \right) \right|\nonumber \\
    &\leq \frac{C_2}{\eta}\left|\beta\eta\mathcal{A}f(0)-\frac{1}{k}\mathcal{A}f(0)\sum_{j=1}^k \int_0^{\varphi_j} v_j du_j - \mathcal{A}f(0)\beta\eta\right|+ o(1)\nonumber \\
    &\leq\frac{1}{k}\left|\mathcal{A}f(0)\right|\sum_{j=1}^k v_j(\varphi_j)+ o(1) \label{eq:22},
\end{align}
where $C_2 = 2C \max_i u_i(r_i)$. The second inequality holds for sufficiently large $n$ (depending on $\eta$). Given that $v_j$ is right-continuous with $v_j(0) = 0$, it follows that $\mathrm{E}_x M\cdot|\mathcal{E}^\infty(\eta)| < \varepsilon/10$ for sufficiently small $\eta$ and large $n$ (depending on $\eta$).

Next, let us show that 
$\mathrm{E}_x M\cdot|\mathcal{E}^n - \mathcal{E}^\infty(\eta)| \rightarrow 0$ as $n \rightarrow \infty$ if $\eta$ is fixed. 
Indeed,

\begin{align}
    &\mathrm{E}_x M\cdot|\mathcal{E}^n - \mathcal{E}^\infty(\eta)| \notag \\
    \leq& \mathrm{E}_x M\cdot\left| 
    \sum_{i=1}^{k} 
    \left[\frac{f(\varphi_i)-f(0)}{u^n_i(\varphi_i)}\right]
    \left(\sum_{j=1}^k \frac{1}{u^n_j(\varphi_j)}\right)^{-1}
    \right. \nonumber \\
    &\left.
    -\sum_{i=1}^{k} 
    \left[\frac{f(\varphi_i)-f(0)}{u_i(\varphi_i)}\right]
    \left(\sum_{j=1}^k \frac{1}{u_j(\varphi_j)}\right)^{-1}
    \right| \nonumber \\ 
    +&\mathrm{E}_x M\cdot|\mathcal{A}f(0)|\cdot
    \left|
    \left(
    \sum_{j=1}^k 
    \frac{1}{u^n_j(\varphi_j)}
    \int_0^{\varphi_j} \bar v_j^n(y)\,du^n_j(y)
    +\sum_{j=1}^k c^n_j+k\beta^n
    \right)
    \left(
    \sum_{j=1}^k \frac{1}{u^n_j(\varphi_j)}
    \right)^{-1}
    \right. \nonumber\\
    &\left.
    -\left(
    \sum_{j=1}^k 
    \frac{1}{u_j(\varphi_j)}
    \int_0^{\varphi_j} v_j(y)\,du_j(y)
    +k\beta
    \right)
    \left(
    \sum_{j=1}^k \frac{1}{u_j(\varphi_j)}
    \right)^{-1}
    \right| \nonumber \\
    =:&\mathrm{E}_x M\cdot
    \left(
    \left|\frac{A^n}{B^n}-\frac{A}{B}\right|
    +
    |\mathcal{A}f(0)|
    \left|\frac{\tilde A^n}{B^n}-\frac{\tilde A}{B}\right|
    \right).
    \label{eq:16}
\end{align} 
Observe that $B >0$ and $A^n \rightarrow A$, $\tilde{A}^n \rightarrow \tilde{A}$, $B_n \rightarrow B$ for each fixed value of $\eta$. Thus, the expression in the parenthesis in \eqref{eq:16} is smaller than $\eta\varepsilon/10C_2$, and then \eqref{eq:16} is bounded by $\varepsilon/10$ by taking a sufficiently large $n$ (depending on $\eta$). The bounds \eqref{eq:19}, \eqref{eq:22} and \eqref{eq:16} together imply that \eqref{2nd-term} is bounded by ${2\varepsilon}/{5}$ for sufficiently small $\eta$ and large $n$. 

Next, let us estimate the term in \eqref{1st-term}. Note that $\eta$ is already fixed.  Let $\Delta_m^n$ denote the term for the inward trip in the $m$-th excursion in \eqref{1st-term}:$$\Delta_m^n := f(X^n_{\tau^n_{m+1} \wedge t}) - f(X^n_{\sigma^n_m \wedge t}) - \int_{\sigma^n_m \wedge t}^{\tau^n_{m+1} \wedge t} \mathcal{A}f(X^n_s) ds.$$

We aim to bound the expectation of the sum of these terms up to the random number of excursions $M$. Using the indicator function $\mathbf{1}_{\{M \geq m\}}$, we rewrite the sum and apply the triangle inequality:$$\left| \mathrm{E}_x \left( \sum_{m=0}^{M} \Delta_m^n \right) \right| 
= \left| \sum_{m=0}^{\infty} \mathrm{E}_x \left( \Delta_m^n \cdot \mathbf{1}_{\{M \geq m\}} \right) \right| 
\leq \sum_{m=0}^{\infty} \left| \mathrm{E}_x \left[ \mathrm{E} ( \Delta_m^n \mid \mathcal{F}_{\sigma_m^n} ) \cdot \mathbf{1}_{\{M \geq m\}} \right] \right|.$$
The term with $m=0$ is 
\[
\left|\mathrm{E}_x\left(f(X^n_{\tau^n_{1} \wedge t})-f(X^n_{0})-\int_{0}^{\tau^n_{1} \wedge t}\mathcal{A}f(X^n_s)ds\right)\right|,
\]
where $\tau^n_1$
is the first arrival to the origin from the initial point. Let us consider the processes $X^n_t$ and $X_t$ on the half-axis $E_i$ that contains the initial point and the origin. Let us extend the processes to the entire real line by continuing the coefficients in such a way that $u^n_i$ converges to $u_i$ and $v^n_i - c^n_i$ converges to $v_i$ in all the points of continuity of $v_i$ on $\mathbb{R}$. One easy way of doing this is by 
setting $u_i^n(x) = u_i(x) = x$ for $x < 0$, while $v_i^n - c_i^n = v_i = x - c$ for $x < 0$, where $c = \sup_n c_i^n \geq 0$ ($c$ is chosen sufficiently large to ensure that all the functions involved are increasing). Momentarily, let us view $X_t^n$ and $X_t$ as processes on the real line rather than the graph.  We can extend $f$ from $E_i$ to the negative half-line in such a way that it belongs to the domain of the generator of $X_t$ (still denoted by $\mathcal{A}$. 
By  \cite{FW}, we have the convergence of $X^n_t$ to  $X_t$.  Since the limiting process is non-degenerate, we have the convergence of the stopped processes as well. Thus,

\begin{align*}
    &\lim_{n \rightarrow \infty} 
\left|\mathrm{E}_x\left(f(X^n_{\tau^n_{1} \wedge t})-f(X^n_{0})-\int_{0}^{\tau^n_{1} \wedge t}\mathcal{A}f(X^n_s)ds\right)\right| \\
=& 
\left|\mathrm{E}_x\left(f(X_{\tau_{1} \wedge t})-f(X_{0})-\int_{0}^{\tau_{1} \wedge t}\mathcal{A}f(X_s)ds\right)\right|= 0.
\end{align*}

That is, for sufficiently large $n$,
\begin{align} \label{tvf}
    &\left|\mathrm{E}_x\left(f(X^n_{\tau^n_{1} \wedge t})-f(X^n_{0})-\int_{0}^{\tau^n_{1} \wedge t}\mathcal{A}f(X^n_s)ds\right)\right|<\frac{\varepsilon}{5}.
\end{align}

When $m>0$, by the strong Markov property at the stopping time $\sigma_m^n$, we can bound the inner expectation:
\begin{align*} \nonumber
\left| \mathrm{E} ( \Delta_m^n \mid \mathcal{F}_{\sigma_m^n} ) \right| \leq \sup_{y \in \Phi} \sup_{t' \leq t} \left| \mathbb{E}_y \left[ f(X^n_{\tau_{1}^n \wedge t'}) - f(y) - \int_{0}^{\tau_{1}^n \wedge t'} \mathcal{A}f(X^n_s)ds \right] \right| =: C(t).
\end{align*}
Substituting this back into the summation, we get:
\begin{align*} 
\left| \mathrm{E}_x \left( \sum_{m=1}^{M} \Delta_m^n \right) \right| 
&\leq C(t) \cdot \sum_{m=1}^{\infty} \mathrm{E}_x \left( \mathbf{1}_{\{M \geq m\}} \right) \nonumber \\
&= C(t) \cdot \sum_{m=1}^{\infty} \mathbb{P}_x(M \geq m)\nonumber \\
&= C(t) \cdot \mathrm{E}_x (M).
\end{align*}
By \eqref{m_boundary}, $\mathrm{E}_x (M) \leq \tilde{c}$ for some constant $\tilde{c}$ independent of $n$ (since $\varphi_i$ and $r_i$ are fixed and $u^n_i$ converge to $u_i$). The quantity $C(t)$ can be bounded by an arbitrarily small constant (e.g., $\varepsilon/(5(\tilde{c}+1)) $) as in \eqref{tvf}. A slight complication in proving the bound on $C(t)$, compared to \eqref{tvf}, is the presence of an extra supremum in $t'$. This, however, is not a problem due to the tightness of the family $X^n_t$.

Thus, \eqref{1st-term} is bounded by
\begin{align} \nonumber
    &\frac{\varepsilon}{5} + 
\frac{\varepsilon}{5(\tilde{c}+1)} \mathrm{E}_x(M)<\frac{2\varepsilon}{5}.
\end{align}

With the bounds on \eqref{1st-term}, \eqref{2nd-term}, and \eqref{3rd-term}, \eqref{inward} is bounded by $\varepsilon$. Therefore, the processes are weakly convergent.

Assume now that the processes converge. The first two conditions, (C1) and (C2), follow directly from  \cite{FW}. Thus we have $u^n_i \to u_i$ and $v_i^n + c_i^n \to v_i$ on $\mathcal{G}\setminus O$, and we only need to verify the third condition (C3).
By weak convergence of the stopped processes, the corresponding exit
times from a fixed neighborhood of the origin converge in distribution. Moreover, such exit times are uniformly integrable.
Therefore, the expectations of the exit times converge, i.e.,    
\[
T^n(0)\longrightarrow T(0),
\]
which translates to

\[\left( \sum_{i=1}^k \frac{1}{u^n_i(\varphi_i)} \int_0^{\varphi_i} v^n_i(y) du^n_i(y) + k\beta^n \right) \left( \sum_{i=1}^k \frac{1}{u^n_i(\varphi_i)} \right)^{-1} \xrightarrow{n\to\infty}  \sum_{i=1}^k  \int_0^{\varphi_i}\frac{1}{k} v_i(y) du_i(y) + \eta\beta,\]
and therefore
$$    \frac{1}{k}\sum_{i=1}^kc_i^n + \beta^n\xrightarrow{n\to\infty}\beta.$$
This completes the proof of the theorem.
\vspace{1cm}
\\
\\
{\bf Acknowledgments}
The authors thank Mark Freidlin for valuable discussions.   Leonid Koralov was supported by the NSF grant DMS-2307377 and the Simons Foundation Grant MP-TSM-00002743.

\end{document}